\documentclass[12pt]{article}

\usepackage{graphicx}
\usepackage{amsfonts}
\usepackage{amssymb,amsmath,amscd}
\usepackage[top=2.75cm, bottom=2.75cm, left=2.75cm, right=2.75cm]{geometry} 

\def\+{{+\!\!\!+}}

\def\d{\partial}

\def\pmb#1{\setbox0=\hbox{#1}%
\kern.0em\copy0\kern-\wd0 
\kern-.04em\copy0\kern-\wd0 
\kern.08em\copy0\kern-\wd0 
\kern-.04em\raise.0433em\box0 }         
 
\newcommand{\nc}{\newcommand} 
\nc{\beq}{\begin{equation}} 
\nc{\eeq}[1]{\label{#1}\end{equation}} 
\nc{\ber}{\begin{eqnarray}} 
\nc{\eer}[1]{\label{#1}\end{eqnarray}} 
\nc{\pek}[1]{\cite{#1}} 
\nc{\enr}[1]{(\ref{#1})} 
\nc{\kal}[1]{{\cal{#1}}} 
\nc{\dott}{\;\cdot\;} 
\nc{\coker}{\mathrm{coker}}
\nc{\ie}{{\it i.e.}}
\nc{\eg}{{\it e.g.}}

\def\0 {\nonumber}

\begin{document} 
\setcounter{page}{0}
\newcommand{\inv}[1]{{#1}^{-1}} 
\renewcommand{\theequation}{\thesection.\arabic{equation}} 
\newcommand{\be}{\begin{equation}} 
\newcommand{\ee}{\end{equation}} 
\newcommand{\bea}{\begin{eqnarray}} 
\newcommand{\eea}{\end{eqnarray}} 
\newcommand{\re}[1]{(\ref{#1})} 
\newcommand{\qv}{\quad ,} 
\newcommand{\qp}{\quad .} 

\thispagestyle{empty}
\begin{flushright} \small
UUITP-06/09  \\
~\\
~\\
{\bf to my wife, Natasha Zabzina}\\
\end{flushright}
\smallskip
\begin{center} \LARGE
{\bf Generalized K\"ahler geometry, \\gerbes, and all that}
 \\[12mm] \normalsize
{\large\bf Maxim Zabzine} \\[8mm]
 {\small\it
 Department of Physics and Astronomy, 
     Uppsala university,\\
     Box 516, 
     SE-75120 Uppsala,
     Sweden
\\~\\
}
\end{center}
\vspace{10mm}
\centerline{\bfseries Abstract} \bigskip
  This work is based on the talk delivered at Poisson 2008. 
 We  review the recent advances in Generalized K\"ahler geometry while  stressing 
  the use of Poisson and symplectic geometry.  The derivation of 
the generalized K\"ahler potential is sketched and the relevant global issues are 
discussed.

\noindent  

\eject
\normalsize



\section{Introduction}
\label{start}

K\"ahler geometry plays a prominent role in mathematical physics. In particular,  it is quite important in modern 
 string theory.  The two dimensional 
 supersymmetric $N=(2,2)$ sigma model should have a  K\"ahler target.  The corresponding quantum theory 
  should be defined over a  Calabi-Yau manifold.  Over the last two decades  the study of these supersymmetric sigma 
   models and their different relatives led to advances in such topics as mirror symmerty,  
   Gromov-Witten invariants and topological strings.
  
 However,  in 1984 it was pointed out by Gates, Hull and Ro\v{c}ek \cite{Gates:1984nk}  that the sigma models 
   with a K\"ahler target are not  the most general supersymmetric $N=(2,2)$ model.  
    They found that the target manifold for these general models should correspond to  bihermitian geometry together with
    some integrability conditions.  The interest in this type of the geometry has been revived after 2002
     when Hitchin introduced the notion of generalized complex structure \cite{Hitchin:2004ut}.  
      In  \cite{Gualtieri:2003dx} Gualtieri gave the alternative description of the Gates-Hull-Ro\v{c}ek geometry
        within the framework of generalized complex geometry and he suggested a new name,  generalized K\"ahler geometry.
     Indeed from the point of view of physics this is a very natural name.   
       There is  hope that many ideas and concepts can be extended to this generalized framework.
      
  In this contribution our goal is modest  and we would like to discuss the different geometrical features of 
    generalized K\"ahler geometry.  We would especially like to stress the Poisson and symplectic aspects of this geometry.
  Our intention will be to review and summarize a number of
    works \cite{Lindstrom:2005zr, Lindstrom:2007qf, Lindstrom:2007xv,  Hull:2008vw}
  written over  a few last years. All of these works were inspired by the tools of supersymmetric sigma models.
   Here we provide the geometrical summary without any reference to sigma models. 

The contribution is organized as follows:   In Section \ref{kahler} we review the standard facts about K\"ahler 
 geometry.  Section \ref{GKG} contains the definition and basic properties of  generalized K\"ahler geometry. 
 In Section \ref{local} we explore the different local description of the geometry and introduce the notion of a generalized
  K\"ahler potential.  Section \ref{global} deals with  ways of gluing the local description and with the interpretation in terms of 
   gerbes.
 Section \ref{end} presents the summary and a list of  open questions. 
 
\section{K\"ahler geometry}
\label{kahler}

Let us remind  the reader of a few well-known facts about  K\"ahler geometry.  In particular we want to discuss the local description of 
 the geometry  and the way of gluing together the local data into a global object.

Consider a complex manifold with a Hermitian metric
$(M, J, g)$. The manifold $M$ is called K\"ahler if the two-form $\omega = g J$ is closed,  $d\omega =0$. The corresponding 
 metric $g$ is called a K\"ahler metric.   The  K\"ahler metrics come in infinite families since on $M$ we can define a new
  closed $2$-form
    $$ \omega' = \omega + i \partial \bar{\partial} \phi~,$$
 which defines another K\"ahler metric $g'$ provided that $\phi$ is a sufficiently small  function. The positivity of the metric 
  is an open condition and thus can be preserved under  small deformations.
  
  Choose an open cover $\{ U_\alpha \}$ of $M$ where all open sets and intersections are contractible. 
 Since $\omega$ is a closed $(1,1)$-form then locally on the patch $U_\alpha$
  we can write 
  \beq
   \omega = i \partial \bar{\partial} K_\alpha~,
  \eeq{kahlerpotential}
   where $K_\alpha (z, \bar{z})$ is a real function on $U_\alpha$ which should gives rise to a positive
    metric. Such a function $K_\alpha$ is called a K\"ahler potential.  Thus locally, provided we choose 
     the complex coordinates $(z, \bar{z})$,  the K\"ahler geometry is defined by any real function which gives 
     rise to a positive metric.
  
    Assume that $\omega/2\pi  \in H^2 (M  , \mathbb{Z})$. The way 
     to glue  the formula (\ref{kahlerpotential})  on the intersection $U_\alpha \cap U_\beta$ is
    \beq
     K_\alpha - K_\beta = F_{\alpha\beta}(z) + \bar{F}_{\alpha\beta} (\bar{z})~,
         \eeq{differencekahpot}
         where $F_{\alpha\beta}(z)$ is a holomorphic function on $U_\alpha \cap U_\beta$.
  Using the fact that $\omega$ is an integral $2$-form we can define the holomorphic 
   transition functions 
   $$G_{\alpha\beta}(z) = e^{F_{\alpha\beta}(z)} ~:~ U_\alpha \cap U_\beta~ \rightarrow~{\mathbb C}_*~,$$
 which satisfy the cocycle condition and the Hermiticity condition         
\beq 
 G_{\alpha\beta} G_{\beta\gamma} G_{\gamma\alpha} =1~,~~~~~~
  G_{\alpha\beta} \bar{G}_{\alpha\beta} = e^{K_\alpha} e^{-K_\beta}~.
\eeq{lineholherm}
    Therefore we are dealing with a  holomorphic line bundle with Hermitian structure. 
     The K\"ahler potential can be defined as
    \beq
     K_\alpha =  \log ||s_\alpha||^2~,
    \eeq{potential}
    where $s_\alpha$ is a local section of a holomorphic line bundle and $||s_\alpha||$ is defined through the Hermitian metric on the line bundle.
     The K\"ahler form $\omega/2\pi$ is the first Chern class of this  holomorphic line bundle with Hermitian structure. 
    This gives us both a local and a global description of the K\"ahler geometry.

\section{Generalized K\"ahler manifolds}
\label{GKG}
 
 Generalized K\"ahler geometry $(M, J_\pm, g, H)$ was introduced originally in \cite{Gates:1984nk}
  as a target manifold for the general $N=(2,2)$ supersymmetric sigma models. The geometry was specified by 
   two complex structures $J_\pm$, a bihermitian metric $g$ and a closed $3$-form $H$ with the following 
    conditions satisfied
 \beq
  \nabla^{\pm} J_\pm =0~,~~~~~~
  \nabla^{\pm} = \nabla \pm g^{-1} H~.
 \eeq{covariancycoc}  
 Equivalently, the generalized K\"ahler geometry can be defined  as 
  a bihermitian manifold $(M, J_\pm, g)$ satisfying the following integrability conditions
 \beq
  d_+^c \omega_+ + d_-^c \omega_ -  = 0~,~~~~~~ d d_{\pm}^c \omega_{\pm}=0~,
 \eeq{Gdjdjdjdj}
 where $\omega_\pm = g J_\pm$ and $d^c = i (\bar{\partial} - \partial)$ with the subscripts "$\pm$"
    referring to the  $J_\pm$ complex structures.   The closed $3$-form $H$ is 
  \beq
   H= d_+^c \omega_+ = - d_-^c \omega_-~.
  \eeq{Hdefkfk}
   The special case $J_+ = J_-$ coincides with the definition of the K\"ahler manifold. 
 The generalized complex description of this bihermitian geometry was given by Gualtieri in
  \cite{Gualtieri:2003dx}.  In the generalized complex language the name "generalized K\"ahler 
   geometry" appears very naturally. In what follows we will not use the language of the generalized 
    geometry, although it appears to be very useful for the discussion of some of the issues. 
 
 The questions we would like to ask are the following: Can we generalize the simple description of 
   K\"ahler geometry reviewed in Section \ref{kahler} to the generalized K\"ahler case? Namely, 
   can we describe the local geometry in terms of a single real function (potential)? If yes, how do we 
    glue them together? In the rest of the contribution we will try to answer these questions.
 
 The definition of generalized K\"ahler geometry can be stated in many different, but equivalent ways. For example, 
  the first condition in (\ref{Gdjdjdjdj})  can be reformulated by saying that the bivectors
  \beq
   \pi_\pm = (J_+ \pm J_-) g^{-1}
 \eeq{realPoisson}
  are Poisson structures  \cite{Lyakhovich:2002kc}.   The Schouten bracket between two 
   Poisson structures defines $H$ as follows
  \beq
   [\pi_+, \pi_-]_s = - 4 g^{-3} H~.
  \eeq{twopoosoll}
  Moreover it has been observed in \cite{Hitchin:2005cv} that the bivector 
  \beq
 \sigma = [J_+, J_-] g^{-1}
 \eeq{definhodllwpp}
  is the real (imaginary) part of the holomorphic Poisson structure with respect to both complex structures. 
   Namely we define $\sigma_\pm = J_\pm \sigma$ to be the imaginary (real) part of these holomorphic Poisson 
    structures. The complex bivector $(\sigma - i \sigma_\pm)$ is a type $(2,0)$  holomorphic bivector  
     for $J_\pm$ complex structure and it is Poisson (Schouten nilpotent). 
      This implies that  $\sigma$ and $\sigma_\pm$ are a pair of real compatible Poisson structures.  Obviously 
       $\sigma$, $\sigma_\pm$ have the same symplectic leaves, although they define different symplectic structures 
        on the leaf.
   The  holomorphic Poisson structures described above are $(2,0)+(0,2)$ parts of the real Poisson structures 
    $\pi_\pm$ \cite{Hitchin:2005cv}.  Thus 
       for the $J_+$ complex structure we have 
 $$ \pi_\pm^{(2,0)} + \pi_\pm^{(0,2)} = \mp \frac{1}{2} J_+ \sigma$$
 and likewise for the  $J_-$ complex structure we have 
 $$ \pi_\pm^{(2,0)} + \pi_\pm^{(0,2)} = \frac{1}{2}  J_-\sigma ~.$$
   Thus we see that there are quite a few Poisson structures on  generalized K\"ahler manifold. Indeed their 
    presence is crucial for the local analysis of the geometry. 
 
\section{Local description}
\label{local}

 In the previous Section we have described two real Poisson structures $\pi_\pm$ and the real part $\sigma$
  of the holomorphic Poisson structure. It is important to stress that $\pi_+$ and $\pi_-$ do not have any common 
   Casimir functions. Moreover the leaf of $\sigma$ is always inside of the leaves for $\pi_\pm$. 
    Indeed the leaves of $\pi_+$ and $\pi_-$ intersect only along a leaf of $\sigma$. 

 Consider a neighborhood of a regular point of a generalized K\"ahler manifold (i.e., there exists a neighborhood of
  the point where the ranks of $\pi_\pm$ are constant).    We can choose 
  the coordinates adapted to the symplectic foliations of the different Poisson structures $\pi_\pm$, $\sigma$
   and complex structures $J_\pm$. Namely we can choose the complex coordinates for $J_+$ 
  \beq
  (z, \bar{z}, z', \bar{z}',  x_+, \bar{x}_+ )~,
  \eeq{coordi+++}
   such that $(x_+, \bar{x}_+ )$ are the coordinates along the leaf of $\sigma$, $(z', \bar{z}',  x_+, \bar{x}_+ )$ are the coordinates 
    along  the leaf of $\pi_-$ and  $(z, \bar{z},  x_+, \bar{x}_+ )$ are the coordinates  along the leaf of $\pi_+$. 
    Analogously we can choose $J_-$ complex coordinates
   \beq
    (z, \bar{z}, z', \bar{z}',  x_-, \bar{x}_-)~,
    \eeq{coordi---}
    such that $(x_-, \bar{x}_-)$ are the coordinates along the leaf of $\sigma$, $(z', \bar{z}',  x_-, \bar{x}_-)$ are the coordinates along the leaf of 
    $\pi_-$ and  $(z, \bar{z},  x_-, \bar{x}_- )$ are the coordinates along the leaf of $\pi_+$.  For these two choices we can pick up the same coordinates 
     along kernels of $\pi_-$ and $\pi_+$.  The possibility of choosing these coordinates follows from the general properties 
      of the Poisson geometry and the definitions (\ref{realPoisson}), (\ref{definhodllwpp}) of Poisson
        structures in terms of the complex structures.  The crucial fact is that these two sets of the coordinates are related 
         to each other by the Poisson diffeomorphism for $\sigma$, i.e. the diffeomorphism preserving $\sigma$. 
  
  \subsection{$\sigma =0$}
  
   We start by considering the special case of generalized K\"ahler geometry when $\sigma =0$ or equivalently, two 
    complex structures commute $[J_+, J_-]=0$.  There exists the integrable local product structure $\Pi = J_+ J_-$ which 
   gives rise to the  real polarization.  We can introduce four differentials: $\d_z$,  $\d_{z'}$ and their complex conjugate 
   $\bar{\d}_{\bar{z}}$, $\bar{\d}_{\bar{z}'}$.  All these differential anticommute with each other. The standard differential 
    we were using before can be written as follows
   $$ d = \d_z + \d_{z'} + \bar{\d}_{\bar{z}} + \bar{\d}_{\bar{z}'}~,~~~~~
   d_+^c = - i\d_z - i \d_{z'}  +i \bar{\d}_{\bar{z}} +i \bar{\d}_{\bar{z}'}~,~~~~~
     d_-^c =  -i \d_z + i\d_{z'} +i \bar{\d}_{\bar{z}} - i \bar{\d}_{\bar{z}'}~.$$
      The corresponding generalized K\"ahler metrics come in infinite families. Namely $2$-forms
   $\omega'_\pm$ on $(M, J_\pm, g)$ 
  \ber
&&   \omega_+' = \omega_+ + i (\d_z \bar{\d}_{\bar{z}} - \d_{z'} \bar{\d}_{\bar{z}'} ) \phi~,\\
&&   \omega_-' = \omega_- + i (\d_z \bar{\d}_{\bar{z}} + \d_{z'} \bar{\d}_{\bar{z}'} ) \phi~,
  \eer{transkkd}
   satisfy the condition (\ref{Gdjdjdjdj}) if the forms $\omega_\pm$ satisfy the same condition.  
    The forms $\omega_\pm'$ define 
    a new bihermitian metric  if $\phi$ is small enough.

  Locally on a patch $U_\alpha$ we can solve the conditions (\ref{Gdjdjdjdj}) as follows
  \beq
 \omega_\pm =  i (\d_z \bar{\d}_{\bar{z}} \mp  \d_{z'} \bar{\d}_{\bar{z}'} )  K_\alpha~,
  \eeq{generalfoe929201}
     where $K_\alpha (z, z', \bar{z}, \bar{z}')$ is a real function such that the corresponding bihermitian 
     metric is positive.  Accordingly,  as result of (\ref{Hdefkfk}) the $3$-form is given 
     \beq
    H = ( \d_z \bar{\d}_{\bar{z}'} \bar{\d}_{\bar{z}}   +
     \d_{z'} \d_z \bar{\d}_{\bar{z}} + \bar{\d}_{\bar{z}} \d_{z'} \bar{\d}_{\bar{z}'}
      + \d_{z'} \d_{z} \bar{\d}_{\bar{z}'}   ) K_\alpha~.
  \eeq{formHbilp}
   This type of generalized K\"ahler geometry is linear generalization of the K\"ahler case. Indeed 
    we are dealing with the local product of two K\"ahler geometries.

  \subsection{invertible $\sigma$}
  
  Now let us consider another special type of generalized K\"ahler geometry when $\sigma$ is invertible.
  Thus $\Omega = \sigma^{-1}$ is a symplectic structure which is the real part of the holomorphic symplectic 
   structure.  Their imaginary parts are given by the corresponding symplectic 
    structures $\Omega_\pm = \Omega J_\pm$. Thus $(\Omega + i \Omega_\pm)$ are the holomorphic 
     symplectic structures  for $J_\pm$.   The symplectic forms $\Omega$, $\Omega_\pm$ 
      encode whole geometry and we can read off from them the complex structures $J_\pm$ and the bihermitian metric. 
  
  The crucial property of a holomorphic symplectic structure is that the complex and Darboux coordinates can be
   chosen simultaneously.  Thus locally we can pick up the Darboux complex coordinates $(q, \bar{q}, p, \bar{p})$
     for $J_+$ such that
  $$  ( \Omega + i \Omega_+) = dq \wedge dp~,$$
   where we choose some polarization. Also we can pick up the Darboux complex 
    coordinates $(Q, \bar{Q}, P, \bar{P})$  for $J_-$ such that
  $$ (\Omega + i \Omega_-) = dQ \wedge dP~,$$
   with some polarization. These two choices of coordinates are related to each other by the symplectomorphism 
    for $\Omega$.  There exist the coordinates $(q, \bar{q}, P, \bar{P})$ and  the 
    generating function $K_\alpha(q, \bar{q}, P, \bar{P})$ such that the corresponding symplectomorphism is 
     defined by the formulas
$$ p = \frac{\d K_\alpha}{\d q}~,~~~~~ \bar{p} = \frac{\d K_\alpha}{\d \bar{q}}~, ~~~~~ 
Q = \frac{\d K_\alpha}{\d P}~, ~~~~~ \bar{Q} = \frac{\d K_\alpha}{\d \bar{P}}~.$$
 Using these expressions we can rewrite the symplectic forms  in the new coordinates $(q, \bar{q}, P, \bar{P})$
  as
\ber
 && \Omega  =  \frac{1}{2}\frac{\d^2 K_\alpha}{\d q \d P} dq \wedge dP + \frac{1}{2}\frac{\d^2 K_\alpha}{\d q \d \bar{P}} dq \wedge d\bar{P} + {\rm c.c.}~,\\
&& \Omega_+  = \frac{i}{2} \frac{\d^2 K_\alpha}{\d q \d \bar{q}} d\bar{q} \wedge d q + 
\frac{i}{2} \frac{\d^2 K_\alpha}{\d \bar{q} \d \bar{P}} d\bar{q} \wedge d\bar{P} + 
  + \frac{i}{2} \frac{\d^2 K_\alpha}{\d \bar{q} \d P} d\bar{q} \wedge d P  - {\rm c.c.}~,\\
 &&\Omega_-  =    \frac{i}{2} \frac{\d^2 K_\alpha}{\d P \d \bar{P}} dP \wedge d \bar{P} +
  \frac{i}{2}  \frac{\d^2 K_\alpha}{\d \bar{P} \d q} dq \wedge d\bar{P} + \frac{i}{2} \frac{\d^2 K_\alpha}{\d \bar{P} \d \bar{q}} d\bar{q} \wedge d\bar{P} 
  - {\rm c.c.}
 \eer{dhdj33oo}
  These are local expressions for $\Omega$, $\Omega_\pm$.  From them we can easily read off  the complex structure 
   and the bihermitian metric. The bihermitian metric can be expressed in terms of the second derivatives of $K$, although 
    the expression is non-linear. Moreover all formulas depend on the choice of polarization $(q, \bar{q}, P, \bar{P})$.
     The polarization can be changed and the generating function should be replaced by the appropriate Legendre transform
      of the original $K_\alpha$.

  \subsection{general case}
  
  The general case can be thought of as a mixture of two previously considered special cases, the linear and 
   non-linear cases.  We will avoid here the full list of  explicit formulas since they are quite lengthy (see \cite{Lindstrom:2005zr} 
   for some of the  explicit expressions).  We just sketch the idea behind their derivation.  As we said before the crucial point is 
    that the complex coordinates for $J_+$ are related to the complex coordinates for $J_-$ through the Poisson 
     diffeomorphism for $\sigma$.  Let the coordinates (\ref{coordi+++}) be 
   \beq
  (z, \bar{z}, z', \bar{z}',  q, \bar{q},  p, \bar{p} )~,
  \eeq{coordi+++djdjj}
   where we choose some polarization along $\sigma$
   and the coordinates (\ref{coordi---})  
   \beq
  (z, \bar{z}, z', \bar{z}',  Q,  \bar{Q}, P, \bar{P})~,
  \eeq{coordi+++dldls}
   with another polarization along $\sigma$.   The coordinates (\ref{coordi+++djdjj}) and (\ref{coordi+++dldls}) are related to each other 
    by the Poisson diffeomorphism for $\sigma$ which can be encoded in the generating function 
      $K_\alpha(z, \bar{z}, z', \bar{z}',  q, \bar{q}, P, \bar{P})$ (as a generating function it has  ambiguities in 
       its definition).   Expressing the complex structures $J_\pm$  in the new coordinates $(z, \bar{z}, z', \bar{z}',  q, \bar{q}, P, \bar{P})$
        through the derivatives of $K_\alpha$
                  one can show that the integrability conditions (\ref{Gdjdjdjdj}) have 
        a solution for $\omega_\pm$ written in terms  of second derivatives of $K_\alpha$.   In the coordinates 
        $(z, \bar{z}, z', \bar{z}',  q, \bar{q}, P, \bar{P})$ the bihermitian metric $g$  can be written in terms of second derivatives of $K_\alpha$. 
         In general the relation will be non-linear in terms of $K_\alpha$. 
   Namely the second relation in (\ref{Gdjdjdjdj})  is solved locally by     
 \beq
  \omega_\pm = d ({\rm Re}\, \lambda_\pm ) + d^c_\pm ({\rm Im}\, \lambda_\pm)~,
 \eeq{definitionlambad}
  where $\lambda_\pm$ are $(1,0)$-forms with respect to $J_\pm$ complex structures.  
   It should be stressed that in (\ref{definitionlambad})  we took into account that $\omega_\pm$ are $(1,1)$-forms
    with respect to the $J_\pm$ complex structures. The first condition in (\ref{Gdjdjdjdj})  implies
  the following compatibility condition between one forms $\lambda_\pm$
 \beq
  d_+^c d ({\rm Re}\, \lambda_+) +   d_-^c d ({\rm Re}\, \lambda_-)=0~. 
 \eeq{integshieiiwow}
  Using the form of the complex structures $J_\pm$ in  the coordinates  $(z, \bar{z}, z', \bar{z}',  q, \bar{q}, P, \bar{P})$
 we can resolve this condition as
 \ber
\label{22dd} {\rm Re}\, \lambda_+& = &\frac{i}{2}( \bar{\d}_{\bar P} + \bar{\d}_{\bar{z}} + \d_{z'}) K_\alpha - {\rm  c.c.}~, \\
 \label{33dd} {\rm Re}\, \lambda_- &=& \frac{i}{2}( \bar{\d}_{\bar q} + \bar{\d}_{\bar{z}} + \bar{\d}_{\bar{z}'}) K_\alpha - {\rm  c.c.} ~,
 \eer{dddd}
  where it is written up to $d$-exact terms which disappear in the final expressions for $\omega_\pm$.  In the expressions 
   (\ref{22dd}) and (\ref{33dd}) we use the locally defined differentials adapted to our coordinates. 
    Now we can read off from (\ref{definitionlambad}), (\ref{22dd}) and (\ref{33dd}) the expression for the bihermitian metric $g$, 
     which will  be in general  non-linear in $K_\alpha$.  The locally defined $2$-forms $d ({\rm Re}\, \lambda_\pm)$ will be 
      non-degenerate if the metric $g$ is non-degenerate. Thus we are dealing with locally defined symplectic structures  $d ({\rm Re}\, \lambda_\pm)$. 
     
     Similar ideas of using Poisson diffeomorphism for $\sigma$ 
      can be utilized in order to generate  new examples of generalized K\"ahler metrics, see
      \cite{hitchindelP}, \cite{gualtieriPois}. 

\section{Global issues vs gerbes}
\label{global}

 In order to understand the global issues we have to figure out how to glue the local formulas  discussed in 
  the previous Section. There are number of complications which we are facing. One of them is the dependence of 
   our formulas on the polarization which we have to pick up on the leaf of $\sigma$ in order to write  everything down.
    The change in the polarization leads to a non-linear Legendre transform of $K_\alpha$ which is unclear how 
     to interpret. The second problem is that we understand only the local description of the generalized K\"ahler geometry 
      in the neighborhood of the regular point and how one deals with the irregular points is unclear to us. 
    
     Below we offer some partial results on the global issues. In the K\"ahler case the holomorphic 
    line  bundles with  Hermitian structure play a central role while in the generalized K\"ahler case the gerbes become important. 
  Gerbes are a geometrical realization of $H^3(M, \mathbb{Z})$ in a manner analogous to the way 
 a line bundle is  geometrical realization of $H^2(M, \mathbb{Z})$ \cite{brylinski}. 
 
 \subsection{biholomorphic gerbe}
 
 The case when $\sigma =0$ is relatively simple one. We have to glue together the local expressions (\ref{generalfoe929201})
  for $\omega_\pm$. On the double intersection $U_\alpha \cap U_\beta$ we have
  \beq
   K_\alpha - K_\beta = f_{\alpha\beta} (z, z') + g_{\alpha\beta}(z, \bar{z'}) + \bar{f}_{\alpha\beta} (\bar{z}, \bar{z}') 
    + \bar{g}_{\alpha\beta}(\bar{z}, z')~,
  \eeq{deje939393}
 where $f_{\alpha\beta} (z, z') $ is $J_+$-holomorphic function on $U_\alpha \cap U_\beta$ and 
   $g_{\alpha\beta}(z, \bar{z'})$  is $J_-$-holomorphic function on $U_\alpha \cap U_\beta$.  
    Assuming that $H \in H^3(M, \mathbb{Z})$ we arrive at the following picture involving the gerbes.
      We can  define over any triple intersections the two sets of transition functions
\beq
G_{\alpha\beta\gamma}(z)~, F_{\alpha\beta\gamma}(z')~:~U_\alpha\cap 
U_\beta \cap U_\gamma~,\rightarrow~\mathbb{C}_*~,
\eeq{definitionbiholgerbe}
which are antisymmetric under permutations of the open sets
and satisfy the cocycle condition on the four-fold intersection. Moreover $G_{\alpha\beta\gamma}(z)$ is 
 holomorphic function with respect to both complex structures, $F_{\alpha\beta\gamma}(z')$ is homolorphic
  for $J_+$ and anti-holomorphic for $J_-$.  We refer to such $G$'s as biholomorphic gerbes and to $F$'s as
   twisted biholomorphic gerbes.  We impose the following "bihermitian" conditions
\beq
G_{\alpha\beta\gamma}  F^{-1}_{\alpha\beta\gamma}= h^+_{\alpha\beta} 
h^+_{\beta\gamma} h^+_{\gamma\alpha}~,~~~~~~~~~~~~~~
G_{\alpha\beta\gamma} \bar{F}_{\alpha\beta\gamma}= 
h^-_{\alpha\beta} h^-_{\beta\gamma} h^-_{\gamma\alpha} ~,
\eeq{bliadi-ia-ustal}
where $h^\pm_{\alpha\beta}$ are $J_\pm$-holomorphic functions 
on double intersections. 
  One can easily see that   the biholomorphic and twisted biholomorphic gerbes 
are are both Hermitian if the conditions (\ref{bliadi-ia-ustal}) are satisfied. 
   From the conditions (\ref{bliadi-ia-ustal}) it follows that there exists  real 
    functions $K_\alpha$ over a patch $U_\alpha$ where
\beq
h^+_{\alpha\beta} \bar{h}^-_{\alpha\beta} 
( h^-_{\alpha\beta})^{-1} ( \bar{h}^+_{\alpha\beta})^{-1} =
e^{K_\alpha} e^{-K_\beta}~.
\eeq{finallypotential}
Comparing with the expression (\ref{deje939393}) we have
 $h_{\alpha\beta}^+ = \exp(f_{\alpha\beta})$ and $h_{\alpha\beta}^- = \exp(g_{\alpha\beta})$. 
  The explicit example of this construction is given by the generalized K\"ahler geometry on
    $S^3 \times S^1$, see \cite{Hull:2008vw}. 

\subsection{general case}

Here we can offer only partial result and some observations.  If we are  dealing with the regular generalized 
 K\"ahler manifold then we can glue the local expressions for $\omega_\pm$ on the double intersections
  $U_\alpha \cap U_\beta$  as follows
 \beq
   K_\alpha - K_\beta = f_{\alpha\beta} (z, z', q) + g_{\alpha\beta}(z, \bar{z'}, P) + \bar{f}_{\alpha\beta} (\bar{z}, \bar{z}', \bar{q}) 
    + \bar{g}_{\alpha\beta}(\bar{z}, z', \bar{P})~,
  \eeq{deje939393}
   where we explicitly ignore the issue of  polarization. Assuming that $H \in H^3(M , \mathbb{Z})$ and proceeding 
    formally we still arrive at the same notion of the biholomorphic gerbe $G_{\alpha\beta\gamma}(z)$ and the twisted biholomorphic 
     gerbe $F_{\alpha\beta\gamma}(z')$ which we have discussed above. 
      The properties (\ref{bliadi-ia-ustal}) and (\ref{finallypotential}) are still satisfied, however 
      now $h_{\alpha\beta}^+(z, z', q) = \exp(f_{\alpha\beta})$ is a $J_+$-holomorphic function of special type (it does not depend on 
       some of the coordinates) 
      and $h_{\alpha\beta}^-(z, \bar{z}', P) = \exp(g_{\alpha\beta})$ is a $J_-$-holomorphic function of special type.

\section{Summary}
\label{end}

Here we presented a discussion of the local and global aspects of the generalized K\"ahler geometry. 
 We reviewed the local description in terms of the generalized K\"ahler potential which is valid in the neighborhood 
  of a regular point. The expression for the bihermitian metric involves the second derivatives of a potential and 
   would be non-linear in general. Thus one can refer to the generalized K\"ahler geometry as  a non-linear 
    generalization of the K\"ahler geometry.  The tools of Poisson geometry  are crucial in the derivations of the present results. 
    
    There are many open  questions which should be addressed.  How to extended the local description to 
     the neighborhood of a irregular point? How  to properly interpret the choice of polarization which is needed
      for the construction to work? In particular it is unclear how to deal with the different choices of the polarization 
       while discussing the global issues.

\bigskip\bigskip
\eject

\noindent{\bf\Large Acknowledgement}:
\bigskip

\noindent I am deeply grateful  to Chris Hull, Ulf Lindstr\"om, Martin Ro\v{c}ek and Rikard von Unge 
 for the collaboration on the present and related subjects.  I thank the organizers for the invitation to
  give this talk.  Finally, I thank the referees for helpful suggestions.
  The research is supported by VR-grant 621-2008-4273.

\end{document}